\documentclass[review]{elsarticle}

\usepackage{lineno,hyperref}
\modulolinenumbers[5]

\journal{Journal of \LaTeX\ Templates}

\newcommand{\be}{\begin{equation}} 
\newcommand{\ee}{\end{equation}}
\newcommand{\beq}{\begin{eqnarray}}
\newcommand{\eeq}{\end{eqnarray}}
\newcommand{\nbeq}{\begin{eqnarray*}}
\newcommand{\neeq}{\end{eqnarray*}}
\newcommand{\D}{\displaystyle}

\begin{document}

\begin{frontmatter}

\title{Characterizations of Exponential Distribution Based on Two-Sided Random Shifts}

\author{Santanu Chakraborty and George P. Yanev }
\address{School of Mathematical and Statistical Sciences \\
   The University of Texas Rio Grande Valley}

\begin{abstract}
A new characterization of the exponential distribution is obtained. 
It is based on an equation involving randomly shifted (translated) order statistics.
No specific distribution is assumed for the shift random variables.
The proof uses a recently developed technique including the Maclaurin series expansion of the probability density of the parent variable.
\end{abstract}


\end{frontmatter}

\section{Introduction and main result}


Let $X_1, X_2, \ldots, X_n$ be a simple random sample from a continuously  distributed parent $X$. Denote by $X_{j:n}$ for $1\le j\le n$ the $j$th order statistic (OS).
Let $n_1$, $n_2$, $k_1$, and $k_2$ be fixed integers, such that $0\le k_i\le n_i-1$ for $i=1,2$. Consider the distributional equation
\be \label{most_general}
X_{n_1-k_1:n_1}+c_1\xi_1 \stackrel{d}{=}X_{n_2-k_2:n_2}+c_2\xi_2,
\ee
where $c_1$ and $c_2$ are certain constants and the "shift" (or "translation")  variables $\xi_1$ and $\xi_2$ are independent from $X_1, X_2, \ldots, X_n$. A variety of particular cases of (\ref{most_general}) have appeared in a number of papers devoted to characterizations of certain classes of continuous distributions. Recent surveys can be found in \cite{AN 16} and \cite{BN 13}.

If $c_2=0$, then (\ref{most_general}) is called one-sided random shift equation. Let $k$ and $n$ be fixed integers, such that $1 \le k\le n-1$. Under some regularity conditions, it is proven in \cite{WA 04} and \cite{CLS 12} that each one of the following two equations  characterizes the exponential distribution:
\[
X_{n-k:n}+\frac{1}{k}\ \xi \stackrel{d}{=} X_{n-k+1:n}\quad \mbox{(consecutive OS)},
\]
\[
X_{n-k:n-1}+\frac{1}{n}\ \xi \stackrel{d}{=} X_{n-k+1:n} \quad
\mbox{(consecutive OS and sample sizes)}
\]
where $\xi$ is unit exponential.

If both $c_1>0$ and $c_2>0$, then we have equation with to two-sided random shifts.
It is established in \cite{WA 04} that for unit exponential $\xi_1$ and $\xi_2$, the equation
\be \label{old_general}
X_{n-k:n-1}+\frac{1}{n}\ \xi_1 \stackrel{d}{=}X_{n-k:n}+\frac{1}{k}\ \xi_2, 
\ee
where $k$ and $n$ are fixed integers, such that $1\le k\le n-1$, characterizes the exponential distribution.

In all equations above the random shifts are  exponentially distributed and therefore (as the characterization implies) are also identically distributed with the parent variable $X$. 
In these notes the condition shifts variables to be exponential is dropped. Instead the weaker assumption that the random shifts are identically distributed with $X$ is made. The one-sided shift case without specifying the distribution of the shift variable is studied in \cite{AV 13} for $n=2$. Namely, it is proven that the equation
\[
X_1+\frac{1}{2}\ X_2 \stackrel{d}{=}X_{2:2}
\]
characterizes the exponential distribution.  This is generalized in \cite{CY 13} and \cite{YC 16} where it is shown that for fixed $n\ge 2$ the equation 
\[
X_{n-1:n-1}+\frac{1}{n}\ X_n \stackrel{d}{=}X_{n:n}, 
\]
characterizes the exponential distribution.

In our main result below we obtain an analog of the two-sided characterization (\ref{old_general}) dropping the assumption on the shifts to be exponential. In the proof of the theorem a recently developed technique based on the Maclaurin series expansion of the density function $f(x)$ is used. This method of proof grew out of an argument given first in \cite{AV 13}.
 
\vspace{0.3cm}{\bf Theorem}\ Let $k$ and $n$ be fixed integers, such that $1\le k\le n-1$. Let $X_1, X_2, \ldots, X_{n+1}$ be a simple random sample from a distribution with cdf $F(x)$ ($F(0)=0$) and pdf $f(x)$ ($f(0)>0$). Assume that $f(x)$ is analytic for  $x>0$. Then for some $\lambda>0$
\be \label{exp}
F(x)=1-e^{-\lambda x}, \qquad x\ge 0
\ee
if and only if 
\be  \label{general}
X_{n-k:n-1}+\frac{1}{n}\ X_n\stackrel{{d}}{=}X_{n-k:n}+\frac{1}{k}\ X_{n+1}.
\ee

In the next section, we present several lemmas needed to prove the theorem. In Section 3, we prove our main theorem.

\section{Preliminaries}

The first lemma plays a key role in the proof of the theorem. Originally it appeared in a similar form as an argument in \cite{AV 13}. In the form presented below, it was proven in \cite{O 15}. For convenience we provide the proof here.

\vspace{0.3cm}{\bf Lemma 1} Let $f(x)$ be analytic for all $x>0$ and $f(0)>0$. If for all non-negative $m$
\be \label{eq109}
f^{(m)}(0)=(-1)^m f^{m+1}(0),
\ee
then
\be \label{exp_pdf}
f(x)=f(0)e^{-f(0)x}.
\ee

{\bf Proof}. It is well-known (cf. \cite{C 78}, p.35) that every (complex) analytic function is infinitely differentiable and, furthermore, has a power series expansion about each point of its domain. Since $f(x)$ is analytic for any $x$, (\ref{eq109}) implies that its Maclaurin series is given  for a positive $x$ by
\[
f(x)=\sum_{j=0}^\infty f^{(j)}(0)\frac{x^j}{j!}=\sum_{j=0}^\infty (-1)^j f^{j+1}(0)\frac{x^j}{j!}=f(0)e^{-f(0)x}.
\]
The proof is complete.

Next, define for all non-negative integers $n$ and $i$, and any real $x$ the numbers
\[
H_{n,i}(x):=\sum_{j=0}^n (-1)^j {n \choose j}(x-j)^i.
\]
It is known (see \cite{R 96}) that $H_{n,n}(\cdot)=n!$ and $H_{n,i}(\cdot)=0$ for $0\leq i \leq n-1$.

In the rest of the section we prove three more lemmas needed in the proof of the theorem.

\vspace{0.3cm}{\bf Lemma 2}\ For $m$, $n$ and $k$ positive integers with $1\le k\le n$,
\be\label{lemma3}
\sum_{l=0}^m {m\choose l} H_{n-k+l,i}(n-k+l+1)
=H_{n-k,i}(n-k+m+1)
\ee
{\bf Proof.}\
Let $s$ and $r$ be positive integers. We shall prove that
\be \label{lemma2}
H_{s,r}(s+1)+H_{s+1,r}(s+2)=H_{s,r}(s+2).
\ee
Using the identity 
\[
{s+1 \choose j} - {s \choose j-1} = {s \choose j}
\]
and the definition of $H_{i,j}(x)$, we obtain for $r\ge s+1$
\nbeq
\lefteqn{\hspace{-1cm}H_{s,r}(s+1)+H_{s+1,r}(s+2)}\\
& & \hspace{-1cm} =
\sum_{j=0}^{s}(-1)^j{s \choose j}(s+1-j)^{r}+\sum_{j=0}^{s+1}(-1)^j{s+1 \choose j}(s+2-j)^{r}\\
& & \hspace{-1cm} =
\sum_{j=1}^{s+1}\left[(-1)^{j-1}{s\choose j-1}+(-1)^j{s+1 \choose j}
\right](s+2-j)^{t}+(s+2)^{r}\\
& & \hspace{-1cm} =
\sum_{j=0}^{s+1}(-1)^{j}{s \choose j}(s+2-j)^{r}\\
& & \hspace{-1cm} = H_{s,r}(s+2).
\neeq
It follows from (\ref{lemma2}) that
\beq 
H_{n-k,i}(n-k+1)+H_{n-k+1,i}(n-k+2)=H_{n-k,i}(n-k+2)\nonumber
\eeq 
Proceeding one more step, we have,
\beq 
H_{n-k,i}(n-k+1)+2H_{n-k+1,i}(n-k+2)+H_{n-k+2,i}(n-k+3)
=H_{n-k,i}(n-k+3)\nonumber
\eeq
Suppose the lemma is true for $m\leq r$. We shall prove it for $m = r+1$. Indeed,
\beq 
\lefteqn{\sum_{l=0}^{r+1} {{r+1}\choose l} 
H_{n-k,l}(n-k+l+1)}\nonumber\\
& & = \sum_{l=0}^r {r\choose l} 
H_{n-k+l,i}(n-k+l+1)+\sum_{l=0}^{r} {r\choose l} 
H_{n-k+l+1,i}(n-k+l+2)\nonumber\\
& & = H_{n-k,i}(n-k+r+1)+H_{n-k+1,i}(n-k+r+2)=H_{n-k,i}(n-k+r+2).\nonumber
\eeq 
The proof is complete.

\vspace{0.3cm} {\bf Lemma 3}\ Let $k$, $n$, and $r$ be positive integers such that $1\le k\le n-1$. Denote $t=n-k-1$. The following two identities are true:
\be \label{sterm}
\sum_{i=0}^r k^i H_{t,t+r+1-i}(n)=\frac{1}{t+1}H_{t+1,t+r+2}(n)-k^{r+1}{t}!
\ee
and
\be \label{fterm} 
\sum_{i=0}^r n^i H_{t,t+r+1-i}(n-1)=\frac{1}{t+1}H_{t+1,t+r+2}(n)-n^{r+1}{t}!. 
\ee

{\bf Proof.}\ 
We shall prove (\ref{sterm}). The proof of (\ref{fterm}) is similar. We have
\beq \label{sterm1}
\lefteqn{\sum_{i=0}^r k^i H_{t,t+r+1-i}(n)
 = \sum_{i=0}^r k^i\sum_{j=0}^{t}(-1)^j {{t}\choose j}(n-j)^{t+r+1-i}}\\
& & = \sum_{j=0}^{t}(-1)^j  {{t}\choose j}\sum_{i=0}^r k^i(n-j)^{t+r+1-i}   \nonumber\\
& & = \sum_{j=0}^{t}(-1)^j  {{t}\choose j}(n-j)^{t+r+1}\sum_{i=0}^r\left(\frac{k}{n-j}\right)^i \nonumber \\
& & = \sum_{j=0}^{t}(-1)^j  {{t}\choose j}(n-j)^{t+r+1}\left[\frac{\D 1-\left(\frac{k}{n-j}\right)^{r+1}}{\D 1-\frac{k}{n-j}}\right]  \nonumber \\
& & = \sum_{j=0}^{t}(-1)^j  {{t}\choose j}\frac{(n-j)^{t+r+2}}{t+1-j}\left[1-\left(\frac{k}{n-j}\right)^{r+1}
    \right]\nonumber \nonumber \\
& & = \frac{1}{t+1}\sum_{j=0}^{t}(-1)^j  {{t+1}\choose j}(n-j)^{t+r+2}\left[1-\left(\frac{k}{n-j}\right)^{r+1}
      \right]\nonumber\\
& & = \frac{1}{t+1}\sum_{j=0}^{t}(-1)^j  {{t+1}\choose j}(n-j)^{t+r+2}-\frac{k^{r+1}}{t+1}\sum_{j=0}^{t}(-1)^j  {{t+1}\choose j}(n-j)^{t+1}  \nonumber \\
& & = \frac{1}{t+1}\sum_{j=0}^{t+1}(-1)^j  {{t+1}\choose j}(n-j)^{t+r+2}-\frac{k^{r+1}}{t+1}\sum_{j=0}^{t+1}(-1)^j  {{t+1}\choose j}(n-j)^{t+1}\nonumber \\
& & =
\frac{1}{t+1}H_{t+1,t+r+2}(n)-\frac{k^{r+1}}{t+1}H_{t+1,t+1}(n)\nonumber \\
& & = \frac{1}{t+1}H_{t+1,t+r+2}(n)-k^{r+1}{t}!. \nonumber
\eeq

\vspace{0.3cm}{\bf Lemma 4} Let $j\ge 1$ and $d$ be integers, such that $j+d\ge 0$. Assume $F(0) = 0$ and for $d\ge 1$ 
\be \label{der}
f^{(m)}(0)=(-1)^m f^{m+1}(0),\qquad m=1,2,\ldots, d.
\ee
Then for $j=1,2,\ldots$
\be \label{eq106}
G_j^{(j+d)}(0)=
\left\{
  \begin{array}{ll}
  H_{j,j+d}(j+1)f^{j+1-d}(0)(f'(0))^d
     & \mbox{if} \quad \quad \qquad d\ge 0;
 \\
    0 & \mbox{if} \quad -j\le d<0.
  \end{array}
\right.
\ee
where  $G_j(x):=F^j(x)f(x)$.

{\bf Proof.}\ 
(i) If $-j\le d<0$, then $G^{(j+d)}_j(0)=0$ because all the terms in the expansion of
$G^{(j+d)}_j(0)$ have a factor $F(0)=0$.

(ii) Let $d=0$. We shall prove (\ref{eq106}) by induction on $j$. One can verify directly the case $j=1$.
Assuming (\ref{eq106}) for $j=k$, we shall prove it for $j=k+1$.
Since $G_{k+1}(x)=F(x)G_{k}(x)$, applying (i), we see that
\nbeq
G^{(k+1)}_{k+1}(0) & = & \sum_{i=0}^{k+1} {k+1 \choose i}F^{(i)}(0) G^{(k+1-i)}_k(0)\\
    & & \hspace{-2cm}=  F(0)G^{(k+1)}_k(0)+(k+1)F'(0)G^{(k)}_k(0)+\sum_{i=2}^{k+1} {k+1 \choose i}F^{(i)}(0) G^{(k+1-i)}_k(0)\\
    & & \hspace{-2cm} = (k+1)!f^{k+2}(0),
\neeq
which completes the proof of (ii).

(iii) Let $d>0$ and $j$ be any positive integer.
 For simplicity, we will write $f^{(i)}:=f^{(i)}(0)$ below.

(a) Let $j=1$. If $d=1$, then we have $G_1^{(2)}(0)=3f'f=f'fH_{1,2}(2)$ since
$H_{1,2}(2)=3$.
Thus, (\ref{eq106}) is true for $d=1$. Next, assuming (\ref{eq106}) for $G_1^{(k)}(0)$, we shall prove it for $G_1^{(k+1)}(0)$.
Since $G_1(x)=F(x)f(x)$, using (\ref{der}) we obtain
\nbeq
G_1^{(k+1)}(0)
    & = & \sum_{i=1}^{k+1} {k+1 \choose i} f^{(i-1)}f^{(k+1-i)}\\
    & = & \sum_{i=1}^{k+1} {k+1 \choose i} (-1)^{i-1}f^i(-1)^{k+1-i}f^{k+2-i}\\
    & = & (-1)f^{k+2}\sum_{j=1}^{k+1}{k+1 \choose j} \\
    & = & (-1)f^{k+2}H_{1,1+k}(2).
    \neeq
This completes the proof for the case (a) $j=1$ and any $d>0$.

(b) Assuming (\ref{eq106}) for $j=1, 2, \ldots k$ and any $d>0$ we shall prove it for $j=k+1$ and any $d>0$.
Since $G_{k+1}(x)=F(x)G_{k}(x)$, by (\ref{der}) and the induction assumption, we obtain
\nbeq
G_{k+1}^{(k+1+d)}(0)
    &  = &
\sum_{i=1}^{k+1+d} {k+1+d \choose i} f^{(i-1)}G_{k}^{(k+1+d-i)}(0)\\
    & = &
  \sum_{i=1}^{d+1} {k+1+d \choose i} f^{(i-1)}G_{k}^{(k+1+d-i)}(0)\\
    & = &
    \sum_{i=1}^{d+1}{k+1+d \choose i}(-1)^{i-1}f^if^{k-d+i}(f')^{1+d-i}H_{k,k+1+d-j}(k+1)\\
    &  = &
    f^{k+2-d}(f')^d \sum_{i=1}^{k+1+d}{k+1+d 
    \choose i}H_{k,k+1+d-i}(k+1)\\
    &  = &
     f^{k+2-d}(f')^d \sum_{l=0}^{k+d}{k+1+d \choose l}H_{k,l}(k+1),
     \neeq 
where in the last equality we have made the index change $l=k+1+d-j$. Therefore, to finish the proof of the induction step (b), we need to show that 
\be \label{bin_H_sum}
\sum_{l=0}^{k+d}{k+1+d \choose l}H_{k,l}(k+1)=H_{k+1,k+1+d}(k+2).
\ee
For brevity denote $r=k+1+d$. Using  the definition of $H_{k,l}(k+1)$, we obtain
\nbeq
\lefteqn{\sum_{l=0}^{r-1}{r \choose l}H_{k,l}(k+1)=\sum_{i=0}^{k}(-1)^i{k \choose i}\sum_{l=0}^{r-1} {r \choose l}(k+1-i)^{l}}\\
& & =\sum_{i=0}^{k}(-1)^i{k \choose i}\left[(k-i+2)^r-(k-i+1)^r\right]\\
    & & =(k+2)^r-\left[(k+1)^r+{k \choose 1}(k+1)^r\right]+ \ldots \\
    & & \hspace{0.5cm} +(-1)^{k}\left[{k \choose k-1}2^r+2^r\right]+(-1)^{k+1}\\
    & & =(k+2)^r - {k+1 \choose 1}(k+1)^r +\ldots +(-1)^{k}{k+1 \choose k}2^r+(-1)^{k+1}\\
    & & =\sum_{j=0}^{k+1} (-1)^j{k+1 \choose j} (k+2-j)^r \\
    & & = H_{k+1,k+1+d}(k+2).
    \neeq
This proves the induction step (b).
Now (iii) follows from (a) and (b). The proof of the lemma is complete.

\section{Proof of the Theorem}
It is not difficult to see that (\ref{general}) is equivalent to (for brevity $t:=n-k-1$)
\[
\int_0^x F^{t}(u)(1-F(u))^{k-1}f(u)f(n(x-u))\, du 
 = \int_0^x F^{t}(u)(1-F(u))^k f(u)f(k(x-u))\, du.
\]
Recalling from Lemma 4 that $G_j(x):=F^j(x)f(x)$ for $j=1,2,\ldots$ and using the binomial formula, we write last  equation as
\nbeq
\lefteqn{\hspace{-2cm}\int_0^x \sum_{l=0}^{k-1}{{k-1}\choose l}(-1)^lG_{t+l}(u)f(n(x-u))\, du}\\
& &  = \int_0^x \sum_{l=0}^k {k \choose l} (-1)^l G_{t+l}(u)f(k(x-u))\, du. \nonumber
\neeq
Differentiating with respect to $x$ as many as $t+r+2$ times for $r\geq 0$ and substituting $x=0$, we obtain
\beq \label{eq23}
\lefteqn{\hspace{-2cm}\sum_{l=0}^{k}(-1)^l{{k-1}\choose l}\sum_{i=0}^{t+r+1}n^i G_{t+l}^{(t+r+1-i)}(0)f^{(i)}(0)}\\
& &  = \sum_{l=0}^k (-1)^l {k \choose l} \sum_{i=0}^{t+r+1} k^i G_{t+l}^{(t+r+1-i)}(0)f^{(i)}(0). \nonumber
\eeq
In view of Lemma 1, to prove that (\ref{eq23}) implies (\ref{exp_pdf}), it is sufficient to show that it implies (\ref{eq109}). We shall prove (\ref{eq109}) by induction with respect to $m$. Let us first verify (\ref{eq109}) for $m=1$, i.e., $f'(0)=-f^2(0)$. If $r=0$ then (\ref{eq23}) becomes
\nbeq 
\lefteqn{\hspace{-1cm}\sum_{l=0}^{k}(-1)^l{{k-1}\choose l}\sum_{i=0}^{t+1}n^i G_{t+l}^{(t+1-i)}(0)f^{(i)}(0)}\\
& &  = \sum_{l=0}^k (-1)^l {k \choose l} \sum_{i=0}^{t+1} k^i G_{t+l}^{(t+1-i)}(0)f^{(i)}(0). \nonumber
\neeq
Since, by Lemma 4, $G_{t+l}^{(t+1-i)}=0$ when $t+1-i<t+l$, omitting the zero terms in the sums and simplifying we obtain
\be \label{eq26}
(t+1)G_{t}^{(t)}(0)f'(0)=-G_{t+1}^{(t+1)}(0)f(0),
\ee
which, in view of Lemma 4, is equivalent to 
\[
(t+1)f^{t+1}(0)t!f'(0)=- f^{t+2}(0)(t+1)!f(0)
\]
and thus $f'(0)=-f^2(0)$. Thus, (\ref{eq109}) is true for $m=1$.

To prove the induction step suppose (\ref{eq109}) holds for $m=1,2,\ldots, r$, i.e.,
\be \label{ind}
f^{(m)}(0)=(-1)^mf^{m+1}(0), \qquad m=1,2,\ldots, r.
\ee
We shall prove it for $m=r+1$. 
 Under the induction hypothesis (\ref{ind}), Lemma~2 implies $G_{t+l}^{(t+r+1-i)}(0)>0$ only if $i\leq r-l+1$ and thus, omitting the zero terms in the sums of (\ref{eq23}), we have

\nbeq 
\lefteqn{\hspace{-2cm}\sum_{l=0}^{k}(-1)^l{{k-1}\choose l}\sum_{i=0}^{r-l+1}n^i G_{t+l}^{(t+r+1-i)}(0)f^{(i)}(0)}\\
& &  = \sum_{l=0}^k (-1)^l {k \choose l} \sum_{i=0}^{r-l+1} k^i G_{t+l}^{(t+r+1-i)}(0)f^{(i)}(0) \nonumber
\neeq
and interchanging the sums, we write it as
\nbeq
\lefteqn{\hspace{-1.7cm}\sum_{i=0}^{r+1}n^i f^{(i)}(0)\sum_{l=0}^{r-i+1}(-1)^l{{k-1}\choose l} G_{t+l}^{(t+r+1-i)}(0)}\\
& &  = \sum_{i=0}^{r+1}k^i f^{(i)}(0)\sum_{l=0}^{r-i+1} (-1)^l {k \choose l}   G_{t+l}^{(t+r+1-i)}(0). \nonumber
\neeq
Collecting in the left-hand side the terms with $i=r+1$, results in 
\beq\label{eq29}
\lefteqn{\hspace{-0.5cm}(n^{r+1}-k^{r+1}) f^{(r+1)}(0)G_{t}^{(t)}(0)}\\
&   = & 
\sum_{i=0}^r f^{(i)}(0)\sum_{l=0}^{r-i+1}\left[k^i{k\choose l}-n^i{k-1 \choose l}\right](-1)^lG_{t+l}^{(t+r+1-i)}(0). \nonumber 
\eeq
Applying Lemma 4 to both sides of (\ref{eq29}) and taking into account the induction  hypothesis (\ref{ind}) in the right-hand side,  it is not difficult to see that
\nbeq
\lefteqn{(n^{r+1}-k^{r+1})f^{(r+1)}(0)t!f^{t+1}(0)}\\
 &  = & (-1)^{r+1}f^{t+r+3}(0)\sum_{i=0}^r\sum_{l=0}^{k}\left[k^i{k\choose l}- n^i{k-1 \choose l}\right]H_{t+l,t+r+1-i}(t+1+l),
\neeq
where $\D {p-1 \choose k}=0$.
Now, it is clear that proving (\ref{eq109}) for $m=r+1$ is equivalent to proving 
\[ \label{eq30}
(n^{r+1}-k^{r+1})t!
=  \sum_{i=0}^r \sum_{l=0}^{k}\left[k^i{k\choose l}
  -n^i{k-1 \choose l}\right]H_{t+l,t+r+1-i}(t+1+l).
\]
Applying  Lemma 2 we write last equation as
\be \label{last_eqn}
(n^{r+1}-k^{r+1})t!
= \sum_{i=0}^r k^i H_{t,t+r+1-i}(n) - \sum_{i=0}^r n^i H_{t,t+r+1-i}(n-1)
\ee
Finally, it is easily to see that (\ref{last_eqn}) with $m=r+1$ follows from (\ref{sterm}) and (\ref{fterm}) in Lemma 3. This proves (\ref{eq109}) for $m=r+1$, which in turn completes the induction argument and theorem's proof.

\end{document}